\documentclass{article}
\usepackage{amsmath}
\usepackage{amssymb}
\usepackage{latexsym}
\usepackage{bbm}
\usepackage{graphicx}
\usepackage{color}
\usepackage{url}
\usepackage[english]{babel}

\newtheorem{thrm}{Theorem}
\newtheorem{prpstn}{Proposition}

\newtheorem{rmrk}{Remark}

\newcommand{\e}{\textrm{e}}

\newcommand{\mbbm}[1]{\mathbbm{#1}}
\newcommand{\ba}{\begin{array}}
\newcommand{\ea}{\end{array}}
\newcommand{\baa}{\left\{ \begin{array}}
\newcommand{\eaa}{\end{array} \right.}
\newcommand{\dt}{\partial_t}
\newcommand{\dx}{\partial_x}
\newcommand{\ds}{\partial_s}

\newcommand{\dxx}{\partial_{xx}}

\newcommand{\dint}{\displaystyle \int}

\newcommand{\dsum}{\displaystyle \sum}
\newcommand{\tld}[1]{\widetilde{#1}}

\newcommand{\Ima}{\textrm{Im\,}}

\newcommand{\Supp}{\textrm{Support }}

\newcommand{\RR}{\mathbb{R}}

\begin{document}
\title{The Back and Forth Nudging algorithm for data assimilation problems: theoretical results on transport equations}
%
\author{Didier Auroux\footnote{Institut de Math\'ematiques de Toulouse, Universit\'e Paul Sabatier Toulouse 3, 31062 Toulouse cedex 9, France; {\sf auroux@math.univ-toulouse.fr}}  \footnote{INRIA, Grenoble, France} \and Ma\"elle Nodet \footnote{Universit\'e de Grenoble, Laboratoire Jean Kuntzmann, UMR 5224, Grenoble, France; {\sf maelle.nodet@inria.fr}} \footnote{INRIA, Grenoble, France} }
\date{\today}

\maketitle
\begin{abstract}
English: \\
In this paper, we consider the back and forth nudging algorithm that has been introduced for data assimilation purposes. It consists of iteratively and alternately solving forward and backward in time the model equation, with a feedback term to the observations.  We consider the case of 1-dimensional transport equations, either viscous or inviscid, linear or not (B\"urgers' equation). Our aim is to prove some theoretical results on the convergence, and convergence properties, of this algorithm. We show that for non viscous equations (both linear transport and B\"urgers), the convergence of the algorithm holds under observability conditions. Convergence can also be proven for viscous linear transport equations under some strong hypothesis, but not for viscous B\"urgers' equation. Moreover, the convergence rate is always exponential in time. We also notice that the forward and backward system of equations is well posed when no nudging term is considered.\\

French:\\
Ce travail \'etudie l'algorithme du nudging direct et r\'etrograde, qui a \'et\'e introduit en assimilation de donn\'ees. Cet algorithme consiste \`a r\'esoudre it\'erativement l'\'equation du mod\`ele, agr\'ement\'ee d'un terme de rappel aux observations, dans le sens direct puis dans le sens r\'etrograde. Dans ce travail nous nous int\'eressons aux \'equations de transport en dimension 1, avec ou sans viscosit\'e, lin\'eaires ou non (B\"urgers). Notre objectif est d'\'etudier la convergence \'eventuelle, et la vitesse de convergence le cas \'ech\'eant, de cet algorithme. Nous prouvons que, pour les \'equations non visqueuses (lin\'eaire ou B\"urgers), la convergence a lieu sous des hypoth\`eses d'observabilit\'e. La convergence peut aussi \^etre d\'emontr\'ee pour des \'equations de transport lin\'eaires visqueuses sous des hypoth\`eses fortes, mais pas pour l'\'equation de B\"urgers visqueuse. En outre, lorsque la convergence a lieu, la vitesse de convergence est toujours exponentielle en temps. Nous remarquons aussi que le syst\`eme d'\'equations directe et r\'etrograde est toujours bien pos\'e lorsqu'aucun terme de rappel n'est pr\'esent.
\end{abstract}
%
%
%
\section{Introduction and main results}
\label{sec:intro}

Data assimilation is the set of techniques aiming to combine in an optimal way the mathematical information provided by the model equations and the physical information given by observations, in order to retrieve the state of a system. Several types of methods have been widely studied in the past decades. We can cite here interpolation, variational and stochastic methods. The first ones interpolate the measurements from the points of observation towards the grid points, the interpolation being weighted by the statistics of the observations \cite{Kalnay}. Variational methods are based on the optimal control theory, and data assimilation is set as being a problem of constrained optimization. The goal is to minimize a cost function measuring the difference between the observations and the corresponding quantities 
 provided by a model integration. The initial condition of the system can then be seen as a control vector \cite{LeDimet}. Finally, the basic idea of stochastic methods is to consider the fields as the realization of a stochastic process and carry out Kalman filtering methods \cite{Kalman,Evensen}. We can also mention one of the very first data assimilation schemes: the nudging method. Also known as Newtonian relaxation or dynamic initialization, it consists of adding a feedback term to the observations directly in the model equations \cite{Anthes}.

All these methods require extensive work, either from the implementation or the computation point of view. For instance, variational methods require the linearization of all operators and also the implementation of the adjoint model. They also need efficient optimization schemes, as the minimization is performed on spaces of huge dimension. On the other side, stochastic methods are somewhat easier to implement, but they require knowledge, storage and manipulations of huge matrices.

The Back and Forth Nudging (BFN) algorithm has recently been introduced as a simple and efficient method for solving data assimilation problems \cite{AurouxBlum1}. In most geophysical applications, data assimilation consists of estimating a trajectory, solution of a partial differential equation (PDE), from the knowledge of observations. These observations are usually sparse in time and space, and incorrect in the sense that they are not the restriction of a solution of the PDE model. One step of the BFN algorithm consists of solving first the model equation, in which a feedback to the observation solution is added, and then the same equation but backwards in time, with also a feedback term to the observations. Such forward and backward integrations provide a new value of the solution at the initial time $t=0$ and the aim of the BFN is to improve the quality of the initial condition.

The idea of the back and forth nudging is to use the difference between the observations and the model trajectory as a feedback control of the equations, both in the forward and backward integrations. This makes the numerical scheme extremely easy to implement, in comparison with both variational and stochastic methods, as we usually only consider diagonal (or even scalar) gain matrices. The back and forth nudging scheme can also be seen as an intermediate scheme between variational and stochastic methods, as the standard nudging technique has both variational (minimization of a compromise between the observations and the energy of the system) and stochastic (sub-optimal Kalman filter) interpretations \cite{AurouxBlum2}.

As a first approximation, we consider in this paper that the observations are correct (i.e. no observation error), and hence the observations satisfy the model equation. We consider various observation domains: first we assume that the observations $u_{obs}(t,x)$ are available for any point $x$ and time $t$, second we assume that they are available for $t\in[t_{1},t_{2}]$ and for all $x$, and third we consider that they are available for all $t$ over a given space domain. This is done through the time and space dependency of the feedback (or nudging) gain matrix $K(t,x)$ that is equal to $0$ when the observations are not available.

Many numerical experiments in almost realistic situations suggest that this algorithm works well, and that the identified solution gets closer to the observations \cite{AurouxBlum2}. The goal of this paper is to prove some theoretical results and convergence properties in the particular case of transport equations, either viscous or inviscid, either linear or non-linear (B\"urgers' equation).

In section \ref{sec:lin-visc}, we consider one step of the BFN algorithm applied to a linear viscous transport equation:
\begin{equation}
\label{eq:1}
\begin{array}{rl}
(F) & \baa{rcl}
\dt u -\nu \dxx u + a(x) \dx u &=& -K (u-u_{obs}) \\
u|_{x=0}=u|_{x=1}&=&0\\
u|_{t=0} &=& u_{0}
\eaa \medskip\\
(B) &
\baa{rcl}
\dt \tld{u} -\nu \dxx \tld{u} + a(x) \dx \tld{u} &=& K^\prime (\tld{u}-u_{obs})\\
\tld{u}|_{x=0}=\tld{u}|_{x=1}&=&0\\
\tld{u}|_{t=T} &=& u(T)
\eaa
\end{array}
\end{equation}
where the following notations hold for all further cases:
\begin{itemize}
\item the time period considered here is $t\in[0,T]$;
\item the first equation $(F)$ is called the forward equation, the second one $(B)$ is called the backward one;
\item $K$ and $K^\prime$ are positive and may depend on $t$ and $x$, but for simplicity reasons, we will always assume that there exists a constant $\kappa \in \RR_{+}^*$ such that $K^\prime(t,x) = \kappa K(t,x)$;
\item $a(x)\in W^{1,\infty}(\Omega)$, $\Omega$ being the considered space domain, either the interval $[0,1]$ or the torus $[0,1]$;
\item $\nu$ is a constant;
\item $u_{obs}$ is a solution of the forward equation with initial condition $u_{obs}^0$:
\begin{equation}
\label{eq:1bis}
\baa{rcl}
\dt u_{obs} -\nu \dxx u_{obs} + a(x) \dx u_{obs} &=& 0 \\
u|_{x=0}=u|_{x=1}&=&0\\
u|_{t=0} &=& u_{obs}^0
\eaa
\end{equation}
\end{itemize}
The following result holds true:
\begin{thrm}
\label{thm:lin-visc} We consider the one-step BFN (\ref{eq:1}) with observations $u_{obs}$ satisfying (\ref{eq:2}). We denote
\begin{equation}
\label{eq:2}
\begin{array}{rcl}
w(t) &=& u(t) - u_{obs}(t) \\
\tld{w}(t) &=& \tld{u}(t)-u_{obs}(t)
\end{array}
\end{equation}
Then we have:
\begin{enumerate}
\item If $K(t,x) = K$, then we have, for all $t\in[0,T]$:
\begin{equation}
\label{eq:3}
\tld{w}(t) = e^{(-K-K^\prime)(T-t)} w(t)
\end{equation}
\item If $K(t,x)=K(x)$, with $\Supp(K) \subset [a,b]$ where $a<b$ and $a\neq 0$ or $b\neq 1$, then equation (\ref{eq:1}) is ill-posed: there does not exist a solution $(u,\tld{u})$, in general.
\item If $K(t,x)=K \mbbm{1}_{[t_{1},t_{2}]}(t)$ with  $0\leq t_{1} < t_{2}\leq T$, then we have
\begin{equation}
\label{eq:4}
\tld{w}(0) = e^{(-K-K^\prime)(t_{2}-t_{1})} w(0)
\end{equation}
\end{enumerate}
\end{thrm}

In section \ref{sec:bg-visc}, we consider one step of the BFN algorithm applied to the viscous B\"urgers' equation:
\begin{equation}
\label{eq:26}
\begin{array}{rl}
(F) & \baa{rcl}
\dt u -\nu \dxx u + u \dx u &=& -K (u-u_{obs}) \\
u|_{x=0}=u|_{x=1}&=&0\\
u|_{t=0} &=& u_{0}
\eaa \medskip\\
(B) &
\baa{rcl}
\dt \tld{u} -\nu \dxx \tld{u} + \tld{u} \dx \tld{u} &=& K^\prime (\tld{u}-u_{obs})\\
\tld{u}|_{x=0}=\tld{u}|_{x=1}&=&0\\
\tld{u}|_{t=T} &=& u(T)
\eaa
\end{array}
\end{equation}
with the same notations as before.\\
The observations $u_{obs}$ satisfy the forward B\"urgers' equation:
\begin{equation}
\label{eq:26bis}
\baa{rcl}
\dt u_{obs} -\nu \dxx u_{obs} + u_{obs} \dx u_{obs} &=& 0 \\
u|_{x=0}=u|_{x=1}&=&0\\
u|_{t=0} &=& u_{obs}^0
\eaa
\end{equation}
We have the following result if $K\ne 0$:
\begin{thrm}
\label{thm:burg-visc}
The BFN iteration (\ref{eq:26}) for viscous B\"urgers' equation, with observations $u_{obs}$ satisfying (\ref{eq:26bis}), is ill-posed, even when $K(t,x)$ is constant (except for $K(t,x)\equiv 0$): there does not exist, in general, a solution $(u,\tld{u})$.
\end{thrm}

In the particular case when $K=K'=0$, the backward problem is ill-posed in the sense of Hadamard, but it has a unique solution if the final condition $\tld{u}|_{t=T}$ is set to a final solution of the direct equation. Moreover, in this particular case, the backward solution is exactly equal to the forward one: $\tld{u}(t) = u(t)$ for all $t\in[0,T]$. The main result is the following:
\begin{prpstn}\label{thrm:K0}
If $K=K'\equiv 0$, then problem (\ref{eq:26}) is well-posed in the sense of Hadamard, and there exists a unique solution $(u,\tld{u})$. Moreover $u=\tld{u}$.
\end{prpstn}

Section \ref{sec:nonvisc} considers the extension of theorem \ref{thm:lin-visc} to the inviscid case, for both linear transport and B\"urgers' equations. \\
We first consider the linear case. The BFN equations are:
\begin{equation}
\label{eq:1nvl}
\begin{array}{rl}
(F) & \baa{rcl}
\dt u + a(x) \dx u &=& -K (u-u_{obs}) \\
u|_{x=0}&=&u|_{x=1}\\
\dx u|_{x=0}&=&\dx u|_{x=1}\\
u|_{t=0} &=& u_{0}
\eaa \medskip\\
(B) &
\baa{rcl}
\dt \tld{u} + a(x) \dx \tld{u} &=& K^\prime (\tld{u}-u_{obs})\\
\tld{u}|_{x=0}&=&\tld{u}|_{x=1}\\
\dx\tld{u}|_{x=0}&=&\dx\tld{u}|_{x=1}\\
\tld{u}|_{t=T} &=& u(T)
\eaa
\end{array}
\end{equation}
where $a(x)$ can be constant or not.

\begin{thrm}
\label{thm:lin-nonvisc} We consider the non viscous one-step BFN (\ref{eq:1nvl}), with observations $u_{obs}$ satisfying (\ref{eq:1nvl}-F) with $K=0$. We denote
\begin{equation}
\label{eq:2nv}
\begin{array}{rcl}
w(t) &=& u(t) - u_{obs}(t) \\
\tld{w}(t) &=& \tld{u}(t)-u_{obs}(t)
\end{array}
\end{equation}
We denote by 
\begin{equation}
\label{mneq:2}
(s,\psi(s,x))
\end{equation}
the characteristic curve of equation (\ref{eq:1nvl}-F) with $K=0$, with foot $x$ in time $s=0$, i.e. such that 
\begin{equation}
\label{mneq:3}
(s,\psi(s,x))|_{s=0} = (0,x)
\end{equation}
We assume that the final time $T$ is such  that the characteristics are well defined and do not intersect over $[0,T]$.\\
Then we have:
\begin{enumerate}
\item If $K(t,x) = K$, 
then we have, for all $t\in[0,T]$:
\begin{equation}
\label{eq:3nvlbis}
\tld{w}(t) = w(t) e^{(-K-K^\prime)(T-t)} 
\end{equation}
\item If $K(t,x)=K \mbbm{1}_{[t_{1},t_{2}]}(t)$ with  $0\leq t_{1} < t_{2}\leq T$, then we have
\begin{equation}
\label{eq:4nvl}
\tld{w}(0) =  w(0) e^{(-K-K^\prime)(t_{2}-t_{1})}
\end{equation}
\item If $K(t,x) = K(x)$, then we have, for all $t\in[0,T]$:
\begin{equation}
\label{eq:3nvl}
\tld{w}(t,\psi(t,x)) = w(t,\psi(t,x))  \,  \exp \left(-\int_{t}^T K(\psi(s,x))+K^\prime(\psi(s,x))\, ds \right)
\end{equation}
\end{enumerate}
\end{thrm}

We finally consider non viscous B\"urgers' equation, still with periodic boundary conditions, and for a time $T$ such that there is no shock in the interval $[0,T]$:
\begin{equation}
\label{eq:1nvb}
\begin{array}{rl}
(F) & \baa{rcl}
\dt u + u \dx u &=& -K (u-u_{obs}) \\
u|_{x=0}&=&u|_{x=1}\\
\dx u|_{x=0}&=&\dx u|_{x=1}\\
u|_{t=0} &=& u_{0}
\eaa \medskip\\
(B) &
\baa{rcl}
\dt \tld{u} + \tld{u} \dx \tld{u} &=& K^\prime (\tld{u}-u_{obs})\\
\tld{u}|_{x=0}&=&\tld{u}|_{x=1}\\
\dx\tld{u}|_{x=0}&=&\dx\tld{u}|_{x=1}\\
\tld{u}|_{t=T} &=& u(T)
\eaa
\end{array}
\end{equation}

\begin{thrm}
\label{thm:bg-nonvisc} We consider the non viscous one-step BFN (\ref{eq:1nvb}), with observations $u_{obs}$ satisfying (\ref{eq:1nvb}-F) with $K=0$. We denote
\begin{equation}
\label{eq:2nvb}
\begin{array}{rcl}
w(t) &=& u(t) - u_{obs}(t) \\
\tld{w}(t) &=& \tld{u}(t)-u_{obs}(t)
\end{array}
\end{equation}
We assume that $u_{obs}\in W^{1,\infty}([0,T]\times\Omega)$, i.e. there exists $M>0$ such that
\begin{equation} \label{eq:2bisnvb}
|\dx u_{obs}(t,x)|\le M, \quad \forall t\in [0,T], \forall x\in \Omega
\end{equation}
Then we have:
\begin{enumerate}
\item If $K(t,x) = K$, then we have, for all $t\in[0,T]$:
\begin{equation}
\label{eq:3nvb}
\|\tld{w}(t)\| \leq e^{(-K-K^\prime+M)(T-t)} \|w(t)\|
\end{equation}
\item If $K(t,x)=K \mbbm{1}_{[t_{1},t_{2}]}(t)$ with  $0\leq t_{1} < t_{2}\leq T$, then we have
\begin{equation}
\label{eq:4nvb}
\|\tld{w}(0)\| \leq e^{(-K-K^\prime)(t_{2}-t_{1})+MT} \|w(0)\|
\end{equation}
\end{enumerate}
\end{thrm}

\begin{prpstn} \label{prpstn:bg-nonvisc}
We consider one forward (resp. backward) BFN step of the non viscous B\"urgers equation (\ref{eq:1nvb}-F) (resp. (\ref{eq:1nvb}-B)). With the notations of theorem \ref{thm:bg-nonvisc}, if $K(t,x)=K(x)$, then we have
\begin{equation} \label{eq:5nvb}
w(T,\psi(T,x)) = w(0,x) \exp \left( -\dint_0^T K(\psi(\sigma,x)) d\sigma-\dint_0^T \partial_x u_{obs}(\sigma,\psi(\sigma,x))d\sigma \right)
\end{equation}
\end{prpstn}

\begin{rmrk} \label{rmk1} For the special case $K(t,x) =K(x) = K \mbbm{1}_{[a,b]}(x)$ where $K$ is a constant and $[a,b]$ is a non-empty sub-interval of $[0,1]$, we have
\begin{equation} \label{eq:5nvbis}
w(T,\psi(T,x)) = w(0,x) \exp \left( -K \chi(x)-\dint_0^T \partial_x u_{obs}(\sigma,\psi(\sigma,x))d\sigma \right)
\end{equation}
where
\begin{equation} \label{eq:6nvb}
\chi(x) = \dint_0^T \mathbbm{1}_{[a,b]}(\psi(\sigma,x)) d\sigma
\end{equation}
is the time during which the characteristic curve $\psi(\sigma,x)$ with foot $x$ of equation (\ref{eq:1nvb}-F) with $K=0$ lies in the the support of $K$. The system is then \emph{observable}  if and only if the function $\chi$ has a non-zero lower bound, i.e. $m := \displaystyle \min_{x} \chi(x) > 0$, the observability being defined by (see \cite{Russell78}):
$$\exists C, \forall u \textrm{ solution of (\ref{eq:1nvb}-F) with } K=0,\quad \|u(T,.)\|^2 \leq C \int_{0}^T \|K(.) u(s,.)\|^2 \, ds
$$
In this case, proposition \ref{prpstn:bg-nonvisc} proves the global exponential decrease of the error, provided $K$ is larger than $\displaystyle \frac{MT}{m}$, where $M$ is defined by equation (\ref{eq:2bisnvb}).
\end{rmrk}

From remark \ref{rmk1}, we can easily deduce that if for each iteration, both in the forward and backward integrations, the observability condition is satisfied, then the algorithm converges. Note that this is not a necessary condition, as even if $\chi(x)=0$, the last exponential of equation (\ref{eq:5nvbis}) is bounded.

Note also that in real geophysical applications (either meteorology or oceanography), there is usually no viscosity. In this case, assuming the observability condition, the BFN algorithm is well posed, and theorem \ref{thm:bg-nonvisc} and proposition \ref{prpstn:bg-nonvisc} say that the solution tends to the observation trajectory everywhere, and not only on the support of $K$. From a numerical point of view, we can observe that even with discrete and sparse observations in space, the numerical solution is corrected everywhere \cite{AurouxBlum2}. We also observed that with a not too large viscosity coefficient, the behavior of the algorithm remains unchanged.

\begin{figure}
\centering
\includegraphics[width=12cm]{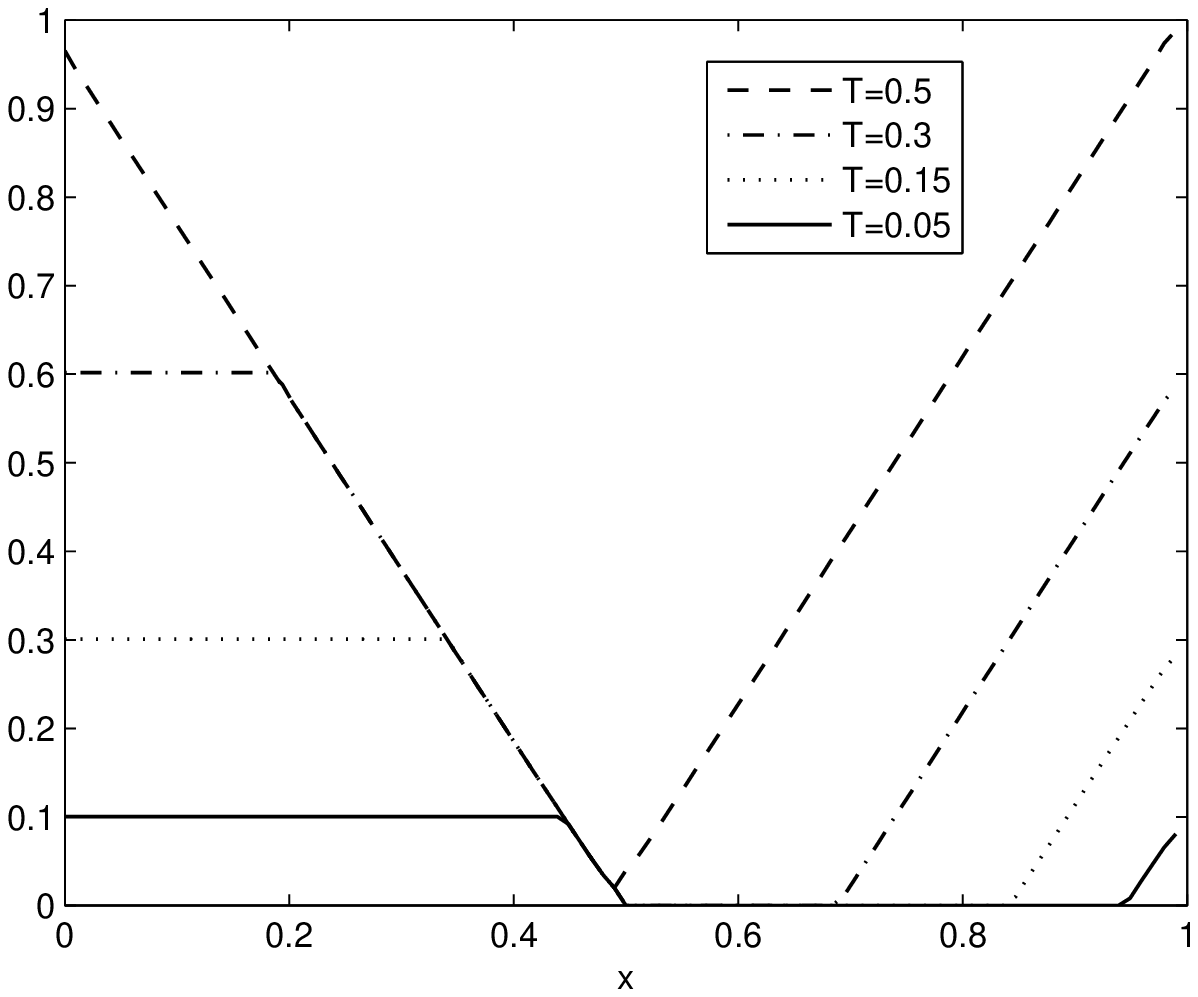}
\includegraphics[width=12cm]{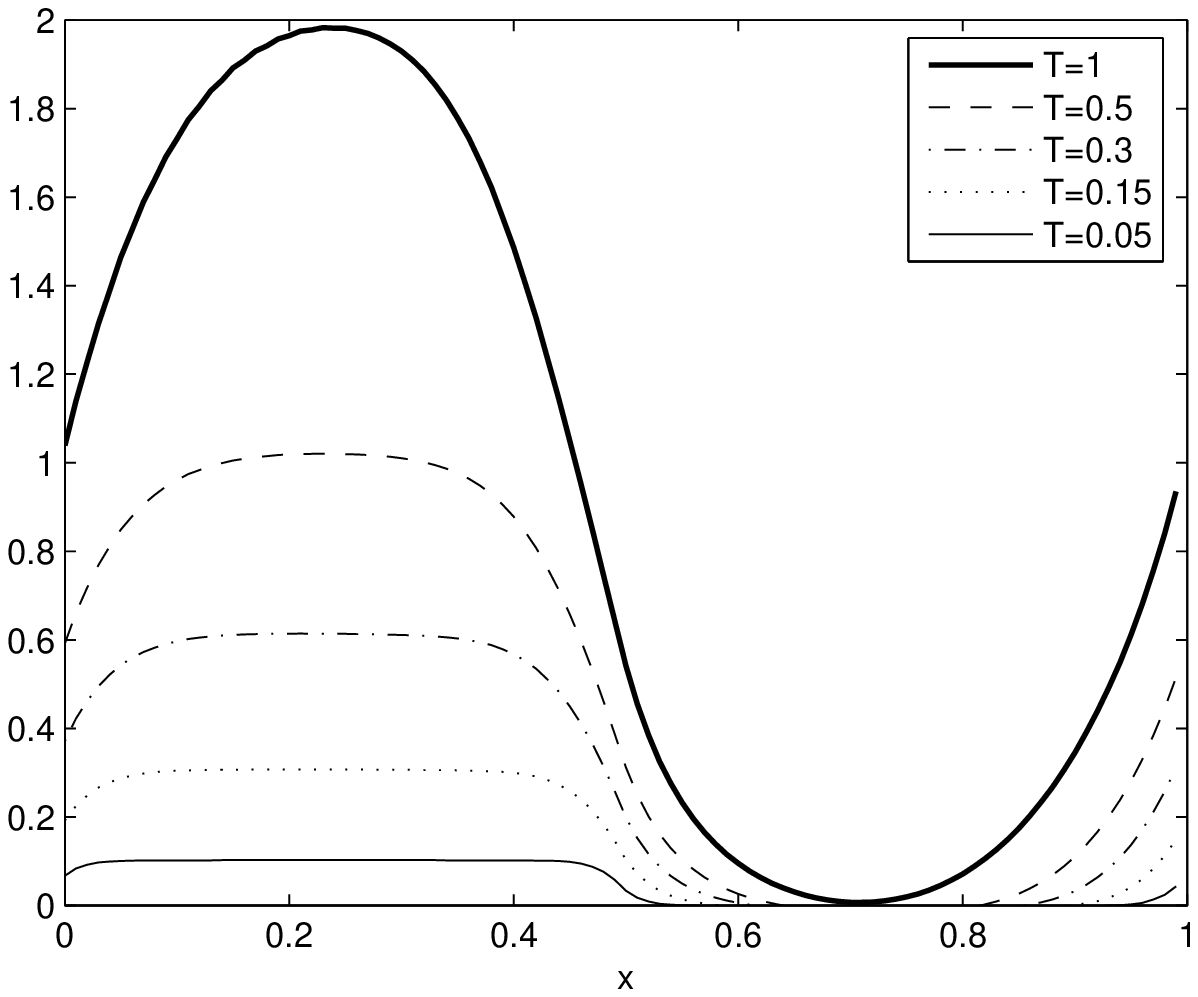}
\caption{Decrease rate of the error after one iteration of BFN (see equation \ref{eq:ww}) as a function of $x$, for various times $T$; top: linear transport equation; bottom: inviscid B\"urgers' equation.}
\label{fig:1}
\end{figure}

Figure \ref{fig:1} illustrates the results given in theorem \ref{thm:lin-nonvisc} in the case $3$ (top) and proposition \ref{prpstn:bg-nonvisc} and remark \ref{rmk1} (bottom). These numerical results correspond to a simple case: $u_{obs}\equiv 0$, $u_0(x) = \alpha \sin(2\pi x)$, $K = K' = \mbbm{1}_{[0;0.5]}(x)$. Various final times $T$ are considered, from $0.05$ to $1$, and both figures show the following expression
\begin{equation}\label{eq:ww}
-\log \left( \frac{\tld{w}(0,x)}{w(0,x)} \right)
\end{equation}
as a function of $x\in [0;1]$. Figure \ref{fig:1}-top illustrates equation (\ref{eq:3nvl}). The best possible decrease rate is then $\max(K+K')\times T=2T$. In the linear case, the transport is $a(x)\equiv 1$. As half of the domain is observed, the observability condition is satisfied iff $T> 0.5$, and this is confirmed by the figure. Concerning B\"urgers' equation, figure \ref{fig:1}-bottom illustrates equation (\ref{eq:5nvb}). After one iteration of BFN, the best possible decrease rate is also $2T$. We can see that in this case, due to the nonlinearities of the model, the solution is less corrected on $[0;0.1]$ but more on $[0.5;0.6]$. From this figure, we can see that the observability condition is satisfied for $T$ larger than approximately $1$.

Finally, some conclusions are given in section \ref{sec:concl}.

\section{Linear transport equation with a viscous term}
\label{sec:lin-visc}

In this section we prove theorem \ref{thm:lin-visc}.

\subsection{Case 1: $K$ constant}
The differences $w$ and $\tld{w}$ satisfy the following equations:
\begin{equation}
\label{eq:5}
\baa{rcl}
\dt w -\nu \dxx w + a(x) \dx w +Kw &=& 0 \\
w|_{x=0}=w|_{x=1}&=&0\\
w|_{t=0} &=& w_{0} \medskip\\
\dt \tld{w} -\nu \dxx \tld{w} + a(x) \dx \tld{w} -K^\prime \tld{w}&=& 0\\
\tld{w}|_{x=0}=\tld{w}|_{x=1}&=&0\\
\tld{w}|_{t=T} &=& w(T)
\eaa
\end{equation}
We denote by $S_{+}$ and $S_{-}$ the operators associated to these equations, seen as forward equations on $[t_{0},t]$ with initial conditions given in $t_{0}$:
\begin{equation}
\label{eq:6}
S_{+}(t_{0},t)(w(t_{0})) = w(t), \quad S_{-}(t_{0},t)(\tld{w}(t_{0})) = \tld{w}(t)
\end{equation}
The BFN algorithm has a solution if and only if we have
\begin{equation}
\label{eq:7}
w(T) \in \Ima(S_{-}(0,T))
\end{equation}
We re-write equation (\ref{eq:6}) associated to $w$:
\begin{equation} 
\label{eq:8}
\baa{rcl}
\dt w -\nu \dxx w + a(x) \dx w -K^\prime w &=& (-K-K^\prime) w \\
w|_{x=0}=w|_{x=1}&=&0\\
w|_{t=0} &=& w_{0}
\eaa
\end{equation}
so that we have, thanks to Duhamel's formula:
\begin{equation}
\label{eq:9}
w(t) = S_{-}(0,t)(w_{0}) + \int_{0}^t S_{-}(s,t)((-K-K^\prime)w(s))\, ds
\end{equation}
If we assume that the expected result is true, i.e. $w(T) \in \Ima(S_{-}(0,T))$, then we can assume that it is also true for all $t$, i.e. we can assume that:
\begin{equation}
\label{eq:10}
\forall t, \exists \varphi(t), w(t) = S_{-}(0,t)\varphi(t)
\end{equation}
In that case, we replace (\ref{eq:10}) in (\ref{eq:9}) and we get:
\begin{equation}
\label{eq:11}
w(t) = S_{-}(0,t)(w_{0}) + \int_{0}^t S_{-}(s,t)((-K-K^\prime)S_{-}(0,s)\varphi(s))\, ds
\end{equation}
As the equation is linear, the scalar coefficient $(K+K^\prime)$ commutes with $S_{-}$ and we get:
\begin{equation}
\label{eq:12}
\begin{array}{rcl}
w(t) &=& S_{-}(0,t)(w_{0}) + (-K-K^\prime)\dint_{0}^t S_{-}(s,t)(S_{-}(0,s)\varphi(s))\, ds\\
&=& S_{-}(0,t)(w_{0}) + (-K-K^\prime)S_{-}(0,t)\dint_{0}^t \varphi(s)\, ds\\
S_{-}(0,t)\varphi(t) &=& S_{-}(0,t)\left[(w_{0}) + (-K-K^\prime)\dint_{0}^t \varphi(s)\, ds\right]
\end{array}
\end{equation}
So that we have
\begin{equation}
\label{eq:13}
\varphi(t) = w_{0}+(-K-K^\prime)\int_{0}^t \varphi(s)\, ds
\end{equation} 
i.e., $\varphi$ satisfies
\begin{equation}
\label{eq:14}
\varphi^\prime(t) = (-K-K^\prime) \varphi, \quad \varphi(0)=w_{0}
\end{equation}
and finally
\begin{equation}
\label{eq:15}
\varphi(t) = w_{0} e^{(-K-K^\prime)t}
\end{equation}
so that we get for $w(T)$:
\begin{equation}
\label{eq:16}
w(T) = S_{-}(0,T)(w_{0}e^{(-K-K^\prime)T})
\end{equation}
Reciprocally, setting 
\begin{equation}
\label{eq:17}
\tld{w}(0) = w_{0} e^{(-K-K^\prime)T}
\end{equation}
leads to $\tld{w}$ satisfying $\tld{w}(T) = w(T)$, so that $(w,\tld{w})$ is the solution of the one-step BFN (\ref{eq:5}).\\
Moreover, we have, for all $t\in [0,T]$:
\begin{equation}
\label{eq:18}
\begin{array}{rcl}
\tld{w}(t) &=& S_{-}(0,t)(\tld{w}(0)) \\
&=& S_{-}(0,t)(w_{0} e^{(-K-K^\prime)T}) \\
&=& e^{(-K-K^\prime)(T-t)}  S_{-}(0,t)(w_{0} e^{(-K-K^\prime)t}) \\
&=& e^{(-K-K^\prime)(T-t)}  S_{-}(0,t)(\varphi(t)) \\
&=& e^{(-K-K^\prime)(T-t)}  w(t)
\end{array}
\end{equation}

\subsection{Case 2: $K(x)$}
We assume that $\Supp(K) \subset [a,b]$ where $a<b$ and $a\neq 0$ or $b\neq 1$, i.e. the support of $K$ is not $[0,1]$. We can follow the same reasoning as previously up to equation (\ref{eq:11}):
\begin{equation}
\label{eq:19}
w(t) = S_{-}(0,t)\varphi(t) =  S_{-}(0,t)(w_{0}) + \int_{0}^t S_{-}(s,t)\left[(-K(x)-K^\prime(x))S_{-}(0,s)\varphi(s)\right]\, ds
\end{equation}
Let us assume, by contradiction, that $-K(x)-K^\prime(x)$ commutes with $S_{-}$. Then we get:
\begin{equation}
\label{eq:20}
S_{-}(0,t)(\varphi(t) - w_{0}) = (-K(x)-K^\prime(x)) S_{-}(0,t)\int_{0}^t  \varphi(s)\, ds
\end{equation}
But we know that $S_{-}$ has the unique continuation property, that is:
\begin{prpstn}
\label{prop:unicity}
If $S_{-}(0,t)(X) = 0$ on a non-empty subset of $[0,1]$, then $S_{-}(0,t)(X) = 0$ on $[0,1]$.
\end{prpstn}
This result and (\ref{eq:20}) give:
\begin{equation}
\label{eq:21}
w(t) = S_{-}(0,t)(\varphi(t)) =  S_{-}(0,t)( w_{0}) = S_{+}(0,t)(w_{0})
\end{equation}
As this stands for every $w_{0}$, we have $S_{-}=S_{+}$ and finally $K=K^\prime=0$, which is a contradiction. Therefore, $K+K^\prime$ does not commute with $S_{-}$. Thus, in general, we cannot find any function $\psi$ such that:
\begin{equation}
\label{eq:22}
\int_{0}^t S_{-}(s,t)\left[(-K(x)-K^\prime(x))S_{-}(0,s)\varphi(s)\right]\, ds = S_{-}(0,t) \psi
\end{equation}

\subsection{Case 3: $K(t)$}

We assume that $K(t,x) = K(t) = K \mbbm{1}_{[t_{1},t_{2}]}(t)$ with $0\leq t_{1}<t_{2}\leq T$. 
We can follow the same reasoning as for $K$ constant, up to the Duhamel formula  (\ref{eq:11}):
\begin{equation}
\label{eq:23}
w(t) = S_{-}(0,t)\varphi(t) =  S_{-}(0,t)(w_{0}) + \int_{0}^t S_{-}(s,t)\left[(-K-K^\prime) \mbbm{1}_{[t_{1},t_{2}]}(s) S_{-}(0,s)\varphi(s)\right]\, ds
\end{equation}
As $K+K^\prime$ is independent of $x$, it commutes with $S_{-}$, and we have for $\varphi$:
\begin{equation}
\label{eq:24}
S_{-}(0,t)\varphi(t) =  S_{-}(0,t)(w_{0}) + 
\baa{ll}
0 & \textrm{if } t \leq t_{1}\smallskip\\
(-K-K^\prime) S_{-}(0,t) \dint_{t_{1}}^t \varphi(s)\, ds  & \textrm{if } t_{1} < t < t_{2}\smallskip\\
(-K-K^\prime) S_{-}(0,t) \dint_{t_{1}}^{t_{2}} \varphi(s)\, ds & \textrm{if }  t \leq t_{2}
\eaa
\end{equation}
So that the corresponding $\varphi$ is given by:
\begin{equation} 
\label{eq:25}
\varphi(t) = 
\baa{ll}
w_{0} & \textrm{if } t \leq t_{1}\smallskip\\
w_{0} e^{(-K-K^\prime)(t-t_{1})} & \textrm{if } t_{1} < t < t_{2}\smallskip\\
w_{0} e^{(-K-K^\prime)(t_{2}-t_{1})}  & \textrm{if }  t \leq t_{2}
\eaa
\end{equation}
And thus the result follows.

\section{B\"urgers' equation with a viscous term}
\label{sec:bg-visc}

\subsection{Proof of theorem \ref{thm:burg-visc}}

Without loss of generality we assume that the observations are identically zero: $u_{obs}(t,x) = 0$ for all $(t,x)$. Let us first introduce some notations.

Let us denote by $w$ (resp. $\tld{w}$) the differences between $u$ (resp. $\tld{u}$) and the observations, as in (\ref{eq:2}), they satisfy the following equations:
\begin{equation}
\label{eq:27}
\begin{array}{rl}
(F) &
\baa{rcl}
\dt w -\nu \dxx w + w \dx w + K w &=& 0\\
w|_{x=0}=w|_{x=1}&=&0\\
w|_{t=0} &=& w_{0} 
\eaa \medskip\\
(B) &
\baa{rcl}
\dt \tld{w} -\nu \dxx \tld{w} + \tld{w} \dx \tld{w} - K^\prime\tld{w} &=& 0\\
\tld{w}|_{x=0}=\tld{w}|_{x=1}&=&0\\
\tld{w}|_{t=T} &=& w(T)
\eaa
\end{array}
\end{equation}
Let us denote also by $S_{+}$ and $S_{-}$ the non-linear operator associated to the forward equations with $K$ or $K^\prime$:
\begin{equation}
\label{eq:28}
S_{+}(t_{0},t)(w(t_{0})) = w(t), \quad S_{-}(t_{0},t)(\tld{w}(t_{0})) = \tld{w}(t), \quad \forall t\geq t_{0}
\end{equation}
We will also use the linear operators $U_{+}$ and $U_{-}$ associated to the following linear equations:
\begin{equation}
\label{eq:29}
\baa{rcl}
\dt \phi  -\nu \dxx \phi + K \phi  &=& 0\\
\phi |_{x=0}=\phi |_{x=1}=0, \qquad
\phi |_{t=0} &=& \phi_{0} 
\eaa
\quad \Longleftrightarrow  \quad U_{+}(0,t)(\phi _{0}) = \phi (t)
\end{equation}
\begin{equation}
\label{eq:29bis}
\baa{rcl}
\dt \phi  -\nu \dxx \phi  -K^\prime \phi  &=& 0\\
\phi |_{x=0}=\phi |_{x=1}=0, \qquad
\phi |_{t=0} &=& \phi_{0} 
\eaa
\quad \Longleftrightarrow  \quad U_{-}(0,t)(\phi _{0}) = \phi (t)
\end{equation}

To prove theorem \ref{thm:burg-visc} we will prove that $w$ is not in the image of $S_{-}$, in general. To do so we will use perturbations theory. We can easily show that $S_{+}$ is infinitely continuous with respect to the data $w_{0}$. So if we suppose that $w_{0}$ is small:
\begin{equation}
\label{eq:30}
w_{0} = \varepsilon \varphi_{0}
\end{equation}
then we have that $w(t)$, solution of the forward equation (\ref{eq:27},$F$) is also small and can be developed in series of $\varepsilon$
\begin{equation}
\label{eq:31}
w = \varepsilon \sum_{n\geq 0} \varepsilon^n w^n
\end{equation}
Similarly, we develop $\tld{w}$ in series of $\varepsilon$
\begin{equation}
\label{eq:31bis}
\tld{w} = \varepsilon \sum_{n\geq 0} \varepsilon^n \tld{w}^n
\end{equation}
As previously, $w$ satisfies:
\begin{equation}
\label{eq:32}
\baa{rcl}
\dt w - \nu \dxx w +K w &=&  - w \dx w\\
w|_{x=0}=w|_{x=1}&=&0\\
w|_{t=0} &=& w_{0} 
\eaa
\end{equation}
so that if we develop in series of $\varepsilon$ we get, for $w^0$: 
\begin{equation}
\label{eq:33}
\baa{rcl}
\dt w^0 - \nu \dxx w^0 +K w^0 &=& 0 \\
w^0|_{x=0}=w^0|_{x=1}&=&0\\
w^0|_{t=0} &=& \varphi_{0} 
\eaa
\end{equation}
For $w^1$ we have:
\begin{equation}
\label{eq:34}
\baa{rcl}
\dt w^1 - \nu \dxx w^1  +K w^1 &=& -w^0 \dx w^0\\
w^1|_{x=0}=w^1|_{x=1}&=&0\\
w^1|_{t=0} &=& 0
\eaa
\end{equation}
Similarly we have for $\tld{w}^0$ and $\tld{w}^1$:
\begin{equation}
\label{eq:33bis}
\baa{rcl}
\dt \tld{w}^0 - \nu \dxx \tld{w}^0  -K^\prime \tld{w}^0 &=& 0 \\
\tld{w}^0|_{x=0}=\tld{w}^0|_{x=1}&=&0\\
\tld{w}^0|_{t=T} &=& w^0(T) 
\eaa
\end{equation}
\begin{equation}
\label{eq:34bis}
\baa{rcl}
\dt \tld{w}^1 - \nu \dxx \tld{w}^1 -K^\prime \tld{w}^1 &=& -\tld{w}^0 \dx \tld{w}^0\\
\tld{w}^1|_{x=0}=\tld{w}^1|_{x=1}&=&0\\
\tld{w}^1|_{t=T} &=& w^1(T)
\eaa
\end{equation}
We can compute $w^0$ and $w^1$ thanks to $U_{+}$:
\begin{equation}
\label{eq:35}
\ba{rcl}
w^0(t) &=& U_{+}(0,t)(\varphi_{0})\\
w^1(t) &=& -\dint_{0}^t U_{+}(s,t)[w^0(s) \dx w^0(s)]\, ds
\ea
\end{equation}
If we assume that $\tld{w}^0$ is well defined, with  
\begin{equation}
\label{eq:35bis}
\tld{w}^0(t)=U_{-}(0,t)(\psi_{0})
\end{equation}
then the condition $\tld{w}(T)=w(T)$ leads to
\begin{equation}
\label{eq:36}
\ba{rcrcl}
&& U_{-}(0,T)(\psi_{0}) &=& U_{+}(0,T)(\varphi_{0})\\
&\Rightarrow& \psi_{0} &=& U_{-}(0,T)^{-1}U_{+}(0,T)(\varphi_{0})\\
&\Rightarrow& \psi_{0} &=& \e^{-(K+K^\prime)T} \varphi_{0}
\ea
\end{equation}
Then we have for $\tld{w}^0$:
\begin{equation}
\label{eq:37}
\ba{rcl}
\tld{w}^0(t) &=& U_{-}(0,t)(\psi_{0})\\
&=& U_{-}(0,t) \e^{-(K+K^\prime)T} \varphi_{0}
\ea
\end{equation}
For $\tld{w}^1$ the final condition $\tld{w}^1(T)=w^1(T)$ gives, thanks to (\ref{eq:35}):
\begin{equation}
\label{eq:38}
\ba{rcl}
\tld{w}^1(T) &=& w^1(T)\\
&=& -\dint_{0}^T U_{+}(s,t)[w^0(s) \dx w^0(s)] \, ds
\ea
\end{equation}
On the other hand, if we assume that $\tld{w}^1$ is well defined, with  $\tld{w}^1(0)=\psi_{T}$, then equation (\ref{eq:34bis}) and the Duhamel formula give
\begin{equation}
\label{eq:39}
\tld{w}^1(T) = U_{-}(0,T)[\psi_{T}]-\int_{0}^T U_{-}(s,T)[\tld{w}^0(s) \dx \tld{w}^0(s)] \, ds\end{equation}
Then, equalling (\ref{eq:39}) and (\ref{eq:38}) we should have
\begin{equation}
\label{eq:40}
U_{-}(0,T)[\psi_{T}]-\int_{0}^T U_{-}(s,T)[\tld{w}^0(s) \dx \tld{w}^0(s)] \, ds \ \ = \ \ -\int_{0}^T U_{+}(s,T)[w^0(s) \dx w^0(s)] \, ds
\end{equation}
Therefore
\begin{equation}
\label{eq:41}
2 U_{-}(0,T)[\psi_{T}]\ \ = \ \ \int_{0}^T U_{-}(s,T)[\dx (\tld{w}^0(s)^2) ] \, ds\   -\int_{0}^T U_{+}(s,T)[\dx (w^0(s) ^2)] \, ds
\end{equation}
If we assume that $\psi_{T} = \displaystyle \frac12 \dx g_{T}$, then we obtain, up to a constant
\begin{equation}
\label{eq:42}
U_{-}(0,T)[g_{T}]\quad = \quad\int_{0}^T U_{-}(s,T)[\tld{w}^0(s)^2] \, ds  -\int_{0}^T U_{+}(s,T)[w^0(s) ^2] \, ds
\end{equation}
We now use (\ref{eq:35}), (\ref{eq:35bis}) and (\ref{eq:36}):
\begin{equation}
\label{eq:43}
\ba{rcl}
U_{-}(0,T)[g_{T}]\quad &=& \quad\dint_{0}^T U_{-}(s,T)[(U_{-}(0,s)(\e^{-(K+K^\prime)T}\varphi_{0})]^2 \, ds  \\
&& -\dint_{0}^T U_{+}(s,T)[(U_{-}(0,s)(\varphi_{0})]^2 \, ds\\
&=& (\e^{-2(K+K^\prime)T}-1) \dint_{0}^T U_{+}(s,T)[(U_{-}(0,s)(\varphi_{0})]^2 \, ds
\ea
\end{equation}
And if $K>0$ and $K^\prime>0$ this last equation is in general impossible: such $g_{T}$ does not, in general, exist. Indeed, let us do an explicit computation thanks to Fourier series:
\begin{equation}
\label{eq:44}
\varphi_{0} = \sum_{n\geq 1}a_{n} \e^{inx}, \quad g_{T} = \sum_{n\geq 1}b_{n} \e^{inx}
\end{equation}
We recall that we have
\begin{equation}
\label{eq:45}
\ba{rcl}
U_{+}(s,t)\left[\dsum_{n\geq 1} c_{n}\e^{inx}\right] &=& \dsum_{n\geq 1} c_{n}\e^{inx} \e^{(-K-\nu n^2)(t-s)}\\
U_{-}(s,t)\left[\dsum_{n\geq 1} c_{n}\e^{inx}\right] &=& \dsum_{n\geq 1} c_{n}\e^{inx} \e^{(K^\prime-\nu n^2)(t-s)}\\
\ea
\end{equation}
Then we can compute the right hand side of equation (\ref{eq:43}):
\begin{equation}
\label{eq:46}
\ba{rcl}
&& (\e^{-2(K+K^\prime)T}-1) \dint_{0}^T U_{+}(s,T)[(U_{-}(0,s)(\varphi_{0})]^2 \, ds\\
&=& (\e^{-2(K+K^\prime)T}-1) \dint_{0}^T U_{+}(s,T) \left[\dsum_{n} a_{n} \e^{K^\prime s} \e^{inx} \e^{-s\nu n^2}\right]^2 \, ds\\
&=& (\e^{-2(K+K^\prime)T}-1) \dint_{0}^T U_{+}(s,T) \left[\dsum_{n} \e^{2sK^\prime } \e^{inx} \sum_{p+q=n}a_{p} a_{q}   \e^{-s\nu (p^2+q^2)}\right] \, ds\\
&=& (\e^{-2(K+K^\prime)T}-1) \dint_{0}^T  \left[\dsum_{n} \e^{-K(T-s)} \e^{2sK^\prime } \e^{-\nu (T-s) n^2}\e^{inx} \sum_{p+q=n}a_{p} a_{q}   \e^{-s\nu (p^2+q^2)}\right] \, ds\\
&=& (\e^{-2(K+K^\prime)T}-1) \dint_{0}^T  \left[\dsum_{n} \dsum_{p+q=n} a_{p} a_{q} \e^{-KT -\nu T n^2 + inx} \e^{2sK^\prime +sK+\nu s n^2-s\nu (p^2+q^2)}\right] \, ds\\
&=& \displaystyle (\e^{-2(K+K^\prime)T}-1)  \left[\dsum_{n} \dsum_{p+q=n} a_{p} a_{q} \e^{-KT -\nu T n^2 + inx} \frac{\e^{2TK^\prime +TK+2\nu pq T } - 1}{2K^\prime+K+2 \nu pq }\right] \\
\ea
\end{equation}
For the left hand side of (\ref{eq:43}) we have:
\begin{equation}
\label{eq:47}
\ba{rcl}
U_{-}(0,T)[g_{T}] &=& \dsum_{n} b_{n}\e^{inx} \e^{K^\prime T -\nu n^2 T}
\ea
\end{equation}
So that we get, for all $n$:
\begin{equation}
\label{eq:48}
\ba{rcl}
b_{n}  &=& \displaystyle \e^{(-K^\prime+K) T}(\e^{-2(K+K^\prime)T}-1)  \left[ \dsum_{p+q=n} a_{p} a_{q}  \frac{\e^{2TK^\prime +TK+2\nu pq T } - 1}{2K^\prime+K+2\nu pq}\right]
\ea
\end{equation}
This defines a distribution  iff $b_{n}$ has polynomial growth, iff $\underline{b_{n}}$ has polynomial growth, where 
\begin{equation}
\label{eq:49}
\ba{rcl}
\underline{b_{n}}  &=& \displaystyle (\e^{-2(K+K^\prime)T}-1)  \left[ \dsum_{p+q=n} a_{p} a_{q}  \frac{\e^{2TK^\prime +TK+2\nu pq T }}{2K^\prime+K+2\nu pq}\right]
\ea
\end{equation}
which is clearly not the case for every sequence $(a_{n})$ with polynomial growth, unless $K=K^\prime=0$.

\subsection{Particular case: $K=K'=0$}

We consider the particular case where $K=K'=0$, i.e. there is no nudging term either in the forward or backward equations. In this case, proposition \ref{thrm:K0} holds true.

Of course, the backward equation itself is ill-posed, as even if there is existence and unicity of the solution (e.g. if the final condition $\tld{u}(T)$ comes from a resolution of the forward equation over the same time period), it does not depend in a continuous way of the data.

The proof is straightforward by using the following Cole-Hopf transformations \cite{Cole,Hopf}:
\begin{equation} \label{eq:colehopf}
\begin{array}{rclcrcl}
u &=& \displaystyle -2\nu \frac{\dx v}{v}, & \quad & v(t,x) = v(t,0) e^{-\frac{1}{2\nu}\dint_0^x u(t,s)\ ds}\\[0.3cm]
\tld{u} &=& \displaystyle -2\nu \frac{\dx \tld{v}}{\tld{v}}, & \quad & \tld{v}(t,x) = \tld{v}(t,0) e^{-\frac{1}{2\nu}\dint_0^x \tld{u}(t,s)\ ds}
\end{array}
\end{equation}
in the forward and backward equations respectively. These transformations allow us to consider the same forward and backward problem, but on the heat equation. The Fourier transform gives the existence and uniqueness of a solution to the forward and backward heat equation, and the equality between the forward $v$ and backward $\tld{v}$ solutions. Equations (\ref{eq:colehopf}) extend the result to the viscous B\"urgers' equation.

\section{Non viscous transport equations}
\label{sec:nonvisc}

\subsection{Linear case: proof of theorem \ref{thm:lin-nonvisc}}

In this section we prove theorem \ref{thm:lin-nonvisc}.

The first two points of the theorem are easily proven as in theorem \ref{thm:lin-visc} with a vanishing viscosity.

Thus we only prove the third point. To do so, we recall that the curves $(s,\psi(s,x))$ are the characteristics of the direct equation (\ref{eq:1nvl}-F) with $K=0$, such that $(s,\psi(s,x))|_{s=0} = (0,x)$ (see \cite{Courant,Evans} for characteristics theory).

For the forward equation  (\ref{eq:1nvl}-F), this change of variable gives
\begin{equation}
\label{mneq:4}
\ds w(s,\psi(s,x)) = -K(\psi(s,x)) w(s,\psi(s,x))
\end{equation}
So that
\begin{equation}
\label{mneq:5}
w(s,\psi(s,x)) =  w(0,x) \, \exp \left(-\int_{0}^s K(\psi(\sigma,x))\, d\sigma \right)
\end{equation}
And in particular for $w(T)$ we have
\begin{equation}
\label{mneq:10}
w(T,\psi(T,x)) =  w(0,x) \, \exp \left(-\dint_{0}^T K(\psi(\sigma,x))\, d\sigma \right)
\end{equation}
For $\tld{w}$ we have similarly
\begin{equation}
\label{mneq:8}
\ds \tld{w}(s,\psi(s,x)) = K^\prime(\psi(s,x)) \tld{w}(s,\psi(s,x))
\end{equation}
So that we have:
\begin{equation}
\label{mneq:9}
\ba{rcl}
\tld{w}(s,\psi(s,x)) &=&  \tld{w}(T,\psi(T,x)) \, \exp \left(-\dint_{s}^T K^\prime(\psi(\sigma,x))\, d\sigma \right)\\
&=&  w(T,\psi(T,x)) \, \exp \left(-\dint_{s}^T K^\prime(\psi(\sigma,x))\, d\sigma \right)
\ea
\end{equation}
Using (\ref{mneq:10}) and (\ref{mneq:5}) we get
\begin{equation}
\label{mneq:11}
\ba{cl}
&\tld{w}(s,\psi(s,x))\\
=&   w(0,x) \, \exp \left(-\dint_{0}^T K(\psi(\sigma,x))\, d\sigma \right)\exp \left(-\dint_{s}^T K^\prime(\psi(\sigma,x))\, d\sigma \right)\\
=&  w(s,\psi(s,x))  \, \exp \left(\dint_{0}^s K(\psi(\sigma,x))\, d\sigma \right)\exp \left(-\dint_{0}^T K(\psi(\sigma,x))\, d\sigma \right) \exp \left(-\dint_{s}^T K^\prime(\psi(\sigma,x))\, d\sigma \right)\\
=&  w(s,\psi(s,x))  \,  \exp \left(-\dint_{s}^T K(\psi(\sigma,x))+K^\prime(\psi(\sigma,x))\, d\sigma \right)
\ea
\end{equation}



\subsection{Non linear case: proof of theorem \ref{thm:bg-nonvisc} and proposition \ref{prpstn:bg-nonvisc}}

From equation (\ref{eq:1nvb}), we deduce that the forward error $w$ satisfies the following equation:
\begin{equation}
\dt w + w\dx w + u_{obs} \dx w + w \dx u_{obs} = -K w
\end{equation}
By multiplying by $w$ and integrating over $\Omega$, we obtain
\begin{equation}
\frac{1}{2}\, \dt \left( \dint_\Omega w^2 \right) + \dint_\Omega w^2 \dx w + \dint_\Omega ( u_{obs} w \dx w + w^2 \dx u_{obs} ) = - \dint_\Omega K w^2
\end{equation}
Some integrations by part give the following:
\begin{equation}
\dt (\| w(t)\|^2) = \dint_\Omega (-2K-\dx u_{obs}) w^2
\end{equation}
We set $M = \| \dx u_{obs} \|_\infty$, and as $K$ does not depend on $x$,
\begin{equation}
\dt (\| w(t)\|^2) \le (-2K+M) \| w(t)\|^2
\end{equation}
We have a similar result for the backward error:
\begin{equation}
\dt (\| \tld{w}(t)\|^2) \le (-2K'+M) \| \tld{w}(t)\|^2
\end{equation}
We first consider the first point of theorem \ref{thm:bg-nonvisc}, i.e. $K(t,x)=K$. Gr\"onwall's lemma between times $t$ and $T$ gives
\begin{eqnarray}
\|w(T)\|^2 &\le& e^{(-2K+M)(T-t)} \|w(t)\|^2 \\
\|\tld{w}(t)\|^2 &\le& e^{(-2K'+M)(T-t)} \|\tld{w}(T)\|^2
\end{eqnarray}
from which equation (\ref{eq:3nvb}) is easily deduced.

In the second case, i.e. $K(t,x) = K \mbbm{1}_{[t_{1},t_{2}]}(t)$ and by successively applying Gr\"onwall's lemma between times $0$ and $t_1$, $t_1$ and $t_2$, and $t_2$ and $T$, one obtains equation (\ref{eq:4nvb}).

Finally, in the case $K(t,x) = K(x)$, by considering a similar approach as in section \ref{sec:nonvisc}.1, i.e. using the characteristics of the direct equation (\ref{eq:1nvb}-F) (resp. B), it is straightforward to prove that
\begin{equation} \label{eq:prp1}
w(s,\psi(s,x)) = w(0,x) e^{-\dint_0^s K(\psi(\sigma,x)) d\sigma} e^{-\dint_0^s \partial_x u_{obs}(\sigma,\psi(\sigma,x)) d\sigma}
\end{equation}
and then,
\begin{equation} \label{eq:prp2}
w(T,\psi(T,x)) = w(0,x) e^{-\dint_0^T K(\psi(\sigma,x)) d\sigma} e^{-\dint_0^T \partial_x u_{obs}(\sigma,\psi(\sigma,x)) d\sigma}
\end{equation}
from which equation (\ref{eq:5nvb}) is easily deduced.

\section{Conclusion}
\label{sec:concl}

Several conclusions can be drawn from all these results. First of all, in many situations, the coupled forward-backward problem is well posed, and the nudging terms allow the solution to be corrected (towards the observation trajectory) everywhere and with an exponential convergence. From a numerical point of view, these results have been observed in several geophysical situations, and many numerical experiments have confirmed the global convergence of the BFN algorithm \cite{AurouxBlum2}.

The second remark is that the worst situation, i.e. for which there is no solution to the BFN problem, is the viscous B\"urgers' equation. But in real geophysical applications, there is most of the time no theoretical viscosity in the equation, and one should consider the inviscid equation instead, for which some convergence results are given. From the numerical point of view, these phenomenon are easily confirmed, as well as the exponential decrease of the error $w$. But we also noticed that if the observations are not too sparse, the algorithm works well even with a quite large viscosity.

Finally, these results extend the theory of linear observers in automatics \cite{Luenberger}: instead of considering an infinite time interval (only one forward equation but for $T\to +\infty$), one can consider an infinite number of BFN iterations on a finite time interval. This is of great interest in almost all real applications, for which it is not possible to consider a very large time period.

\subsection*{Acknowledgement}
The authors are thankful to Prof. G. Lebeau (University of Nice Sophia-Antipolis) for his fruitful ideas and comments. This work has been partially supported by ANR JCJC07 and INSU-CNRS LEFE projects.

\bibliographystyle{alpha}
\bibliography{BFN-burgers-thms}

\end{document}